\documentclass{article}
\usepackage[dvips]{graphicx}
\usepackage[latin1]{inputenc}
\usepackage{amssymb,amsmath,array}
\newtheorem{Lem}{\underline {Lemma}}
\newtheorem{Prop}{\underline {Proposition}}
\newtheorem{Theorem}{\underline {Theorem}}
\def\keywords#1{\gdef\@motcle{#1}}\gdef\@motcle{}

\begin{document}

\title{Efficient Estimation of Multidimensional Regression Model using Multilayer Perceptrons}

%***********************************************************************
% AUTHORS INFORMATION AREA
%***********************************************************************
\author{Joseph Rynkiewicz $^1$
\vspace{.3cm}\\
%
% Addresses and institutions (remove "1- " in case of a single institution)
Université Paris I - SAMOS/MATISSE\\
90 rue de Tolbiac, Paris - France\\ 
}
%***********************************************************************
% END OF AUTHORS INFORMATION AREA
%***********************************************************************

\maketitle

\begin{abstract}
This work concerns the estimation of multidimensional nonlinear regression
models using multilayer perceptrons (MLPs). The main problem with such models is that we need to know the covariance matrix of the noise to get an optimal estimator. However, we show in this paper that if we choose as the cost function  the logarithm of the determinant of the empirical error covariance matrix, then we get an asymptotically optimal estimator. Moreover, under suitable assumptions, we show that this cost function leads to a very simple asymptotic law for testing the number of parameters of an identifiable MLP. Numerical experiments confirm the theoretical results. 
\end{abstract}
\begin{paragraph}{\bf keywords}
non-linear regression, multivariate regression, multilayer Perceptrons, asymptotic normality 
\end{paragraph}
\section{Introduction}
Let us consider a sequence $\left(Y_{t},Z_{t}\right)_{t\in \mathbb{N}}$
of i.i.d. (i.e. independent, identically distributed) random vectors, with $Y_t$ a d-dimensional vector. Each couple $\left(Y_{t},Z_{t}\right)$ has the same law as a generic variable $(Y,Z)$, but it is not hard to generalize all that we show in this paper for stationary mixing variables and therefore for time series. We assume that the model can be written as
\[
Y_{t}=F_{W^0}(Z_{t})+\varepsilon _{t}
\]
where 

\begin{itemize}
\item $F_{W^0}$ is a function represented by an MLP with parameters or
weights $W^0$.
\item $(\varepsilon _{t})$ is an i.i.d.-centered noise with unknown invertible
covariance matrix $\Gamma _{0}$.
\end{itemize}
This corresponds to multivariate non-linear least square model, as in chapters 3.1 and 5.1 of Gallant \cite{Gallant}. Indeed, an MLP function can be seen as a parametric non-linear function, for example an one hidden layer MLP using hyperbolic tangent as transfert functions ($\tanh$) can be written $F_{W^0}(Z_{t})=\left(F^1_{W^0}(Z_{t}),\cdots,F^d_{W^0}(Z_{t})\right)^T$,  where  $T$ denotes the transposition of the matrix,  with :
\[
F^i_{W^0}(z)\sum_{j=1}^Ha_{ij}\tanh\left(\sum_{k=1}^Lw_{jk}z_k+w_{j0}\right)+a_{i0}
\] 
where $H$ is the number of hidden units and $L$ is the dimension of the input $z$, then the parameter vector is 
\[
\left(a_{10},\cdots,a_{dH},w_{10},\cdots,w_{HL}\right)\in{\mathbb R}^{(H+1)\times d+(L+1)\times H}
\]

There are some obvious transformations that can be applied to an MLP without changing its input-output map. For instance, suppose we pick an hidden node $j$ and we change the sign of all the weights $w_{ij}$ for $i=0,\cdots,H$, and also the sign of all $a_{ij}$ for $i=0,\cdots,d$. Since $\tanh$ is odd, this will not alter the contribution of this node to the total net output. Another possibility is to interchange two hidden nodes, that is, to take two hidden nodes $j_1$ and $j_2$ and relabel $j_1$ as $j_2$ and $j_2$ as $j_1$, taking care to also relabel the corresponding weights.  These transformations form a finite group (see Sussmann \cite{Sussmann}). 

We will consider equivalence classes of one hidden layer MLPs: two MLPs are in the same class if the first one is the image by such transformation of the second one, the considered set of parameters is then the quotient space of parameters by this finite group.  In this space, we assume that the model is identifiable it means that the true model belongs to the considered family of models and that we consider MLPs without redundant units. This is a very strong assumption but it is known that estimated weights of an MLPs with redundant units can have a very strange asymptotic behavior (see Kukumizu \cite{Fukumizu2}), because the Hessian matrix is singular. The consequence of the identifiability of the model is that the Hessian matrix computed in the sequel will be definite positive (see Fukumizu \cite{Fukumizu}).  In the sequel we will always assume that we are under the assumptions making the Hessian matrix definite positive. 

\subsection{Efficient estimation}
A popular choice for the associated cost function is the mean square error: 
\begin{equation}\label{ols}
\frac{1}{n}\sum _{t=1}^{n}\left\Vert Y_{t}-F_{W}\left(Z_{t}\right)\right\Vert ^{2}
\end{equation}
 where $\left\Vert .\right\Vert $ denotes the Euclidean norm on $\mathbb{R}^{d}$.
Although this function is widely used, it is easy to show that we then get a suboptimal estimator, with a larger asymptotic variance that the estimator minimizing the generalized mean square error :
\begin{equation}\label{true}
\frac{1}{n}\sum _{t=1}^{n}\left(Y_{t}-F_{W}\left(Z_{t}\right)\right)^T\Gamma_0^{-1}\left(Y_{t}-F_{W}\left(Z_{t}\right)\right)
\end{equation}
But, we need to know the true covariance matrix of the noise to use this cost function.
A possible solution is to use an approximation $\Gamma$ of the covariance error matrix $\Gamma_0$ to compute the generalized least squares estimator :
\begin{equation}
\frac{1}{n}\sum _{t=1}^{n}\left(Y_{t}-F_{W}\left(Z_{t}\right)\right)^{T}\Gamma ^{-1}\left(Y_{t}-F_{W}\left(Z_{t}\right)\right)
\end{equation}
A  way to construct a sequence of $\left(\Gamma_{k}\right)_{k\in \mathbb{N}^{*}}$ yielding
a good approximation of $\Gamma _{0}$ is the following: using the
ordinary least squares estimator $\hat{W}^1_{n}$, the noise covariance
can be approximated by 
\begin{equation}
\Gamma _{1}:=\Gamma \left(\hat{W}^1_{n}\right):=\frac{1}{n}\sum _{t=1}^{n}(Y_{t}-F_{\hat{W}^1_{n}}(Z_{t}))(Y_{t}-F_{\hat{W}^1_{n}}(Z_{t}))^{T}.
\end{equation}
then, we can use this new covariance matrix to find a generalized
least squares estimator $\hat{W}_{n}^{2}$: 
\begin{equation}
\hat{W}_{n}^{2}=\arg \min _{W}\frac{1}{n}\sum _{t=1}^{n}\left(Y_{t}-F_{W}\left(Z_{t}\right)\right)^{T}\left(\Gamma _{1}\right)^{-1}\left(Y_{t}-F_{W}\left(Z_{t}\right)\right)
\end{equation}
and calculate again a new covariance matrix \[
\Gamma _{2}:=\Gamma \left(\hat{W}_{n}^{2}\right)=\frac{1}{n}\sum _{t=1}^{n}(Y_{t}-F_{\hat{W}_{n}^{2}}(Z_{t}))(Y_{t}-F_{\hat{W}_{n}^{2}}(Z_{t}))^{T}.\]
It can be shown  that this procedure gives a sequence of parameters 
\[
\hat{W}_{n}\rightarrow \Gamma _{1}\rightarrow \hat{W}_{n}^{2}\rightarrow \Gamma _{2}\rightarrow \cdots
\]
 minimizing the logarithm of the determinant of the empirical covariance
matrix (see chapter 5 in Gallant\cite{Gallant}) :
 
\begin{equation}\label{logdet}
U_{n}\left(W\right):=\log \det \left(\frac{1}{n}\sum _{t=1}^{n}(Y_{t}-F_{W}(Z_{t}))(Y_{t}-F_{W}(Z_{t}))^{T}\right)
\end{equation}
The use of this cost function for neural networks has been introduced by Williams in 1996 \cite{Williams}, however its theoretical and practical properties have not yet been studied. 
Here, the calculation of the asymptotic properties of $U_{n}\left(W\right)$ will show that this cost function leads to an asymptotically optimal estimator, with the same asymptotic variance that the estimator minimizing (\ref{true}), we say then that the estimator is ``efficient''.

\subsection{testing the number of parameters}
Let $q$ be an integer less than $s$, we want to test ``$H_0: W\in \Theta_q \subset \mathbb R^q$'' against  ``$H_1: W\in \Theta_s \subset \mathbb R^s$'', where the sets $\Theta_q$ and $\Theta_s$ are compact and $\Theta_q\subset\Theta_s$. $H_0$ expresses the fact that $W$ belongs to a subset $\Theta_q$ of $\Theta_s$ with a parametric dimension lesser than $s$ or, equivalently, that $s-q$ weights of the MLP in $\Theta_s$ are null. 
If we consider the classical mean square error cost function: $V_n(W)=\sum_{t=1}^n\Vert Y_{t}-F_{W}(Z_{t})\Vert^2$, we get the following test statistic: 
\[
S_n=n\times\left(\min_{W\in \Theta_q}V_n(W)-\min_{W\in \Theta_s}V_n(W)\right)
\]
Under the null hypothesis $H_0$, it is shown in Yao \cite{Yao} that  $S_n$ converges in law to a weighted sum of $\chi^2_1$
\[
S_n\stackrel{\cal D}{\rightarrow}\sum_{i=1}^{s-q}\lambda_i\chi_{i,1}^2
\]
where the $\chi_{i,1}^2$ are $s-q$ i.i.d. $\chi^2_1$ variables and $\lambda_i$ are strictly positives eigenvalues of the asymptotic covariance matrix of the estimated weights, different from 1 if the true covariance matrix of the noise is not the identity matrix. So, in the general case, where the true covariance matrix of the noise is not the identity matrix, the asymptotic distribution is not known, because the $\lambda_i$s are not known and it is difficult to compute the asymptotic level of the test. 

However, if we use the cost function
\(
U_{n}\left(W\right)
\)
then, under $H_0$, the test statistic:
\begin{equation}\label{tn}
T_n=n\times\left(\min_{W\in \Theta_q}U_n(W)-\min_{W\in \Theta_s}U_n(W)\right)
\end{equation}
will converge to a classical $\chi^2_{s-q}$ so the asymptotic level of the test will be very easy to compute. This is another advantage of using the cost function in Eq. (\ref{logdet}). Note that this result  is true even if the noise is not Gaussian (it is more general that the maximum likelihood estimator) and without knowing the true covariance of the noise $\Gamma_0$, so without using the cost function (\ref{true}) or even an approximation of it.

In order to prove these properties, the paper is organized as follows. First we compute the first and second derivatives of $U_n(W)$ with respect to the weights of the MLP, then we deduce the announced properties with classical statistical arguments. Finally, we confirm the theoretical results with numerical experiments. 

\section{The first and second derivatives of $W\longmapsto U_{n}\left(W\right)$}
First, we introduce a notation:  if $F_W(X)$ is a $d$-dimensional parametric function depending on a parameter vector $W$, let us write $\frac{\partial F_W(X)}{\partial W_k}$ (resp. $\frac{\partial^2 F_W(X)}{\partial W_k\partial W_l}$) for the $d$-dimensional vector of partial derivatives (resp. second order partial derivatives) of each component of $F_W(X)$. Moreover, if $\Gamma(W)$ is a matrix depending on $W$, let us write $\frac{\partial}{\partial W_k}\Gamma(W)$ the matrix of partial derivatives of each component of $\Gamma(W)$. 
\subsection{First derivatives}
Now, if $\Gamma_n(W)$ is a matrix depending on the parameter vector $W$, we get  (see Magnus and Neudecker \cite{Magnus})
\[
\frac{\partial }{\partial W_{k}}\log \det \left(\Gamma _{n}(W)\right)=tr\left(\Gamma_{n}^{-1}(W) \frac{\partial }{\partial W_{k}}\Gamma _{n}(W)\right).\]
Here 
\[
\Gamma_n(W)=\frac{1}{n}\sum _{t=1}^{n}(y_{t}-F_{W}(z_{t}))(y_{t}-F_{W}(z_{t}))^{T}.
\] 
Note that this matrix $\Gamma_{n}(W)$ and it inverse  are symmetric. 
Now, if we note that
\[
A_n(W_k)=\frac{1}{n}\sum _{t=1}^{n}\left(-\frac{\partial F_{W}(z_{t})}{\partial W_{k}}(y_t-F_{W}(z_{t}))^T\right),
\]
then, using the fact 
\[
tr\left(\Gamma_{n}^{-1}(W)A_n(W_k)\right)=tr\left(A_n^T(W_k)\Gamma_{n}^{-1}(W)\right)=tr\left(\Gamma_{n}^{-1}(W)A_n^T(W_k)\right),
\]
we get 
\begin{equation}\label{first_deriv}
\frac{\partial }{\partial W_{k}}\log \det \left(\Gamma _{n}(W)\right)=2tr\left(\Gamma_{n}^{-1}(W)A_n(W_k)\right).
\end{equation}

\subsection{Calculus of the derivative of $W\longmapsto U_{n}\left(W\right)$ for an MLP}
Let us note $\left(\Gamma _{n}\left(W\right)\right)_{ij}$ (resp. $\left(\Gamma^{-1}_{n}\left(W\right)\right)_{ij}$) the element of the $i$th line and $j$th column of the matrix $\Gamma _{n}\left(W\right)$ (resp. $\Gamma^{-1}_{n}\left(W\right))$. We note also $F_{W}(z_{t})(i)$ the $i$th component of a multidimensional function and for a matrix $A=\left(A_{ij}\right)$, we note that $\left(A_{ij}\right)_{1\leq i,j\leq d}$ is the vector obtained by concatenation of the columns of $A$. Following the previous results, we can write for the derivative of $\log (\det (\Gamma _{n}\left(W\right)))$
with respect to the weight $W_{k}$: \[
\frac{\partial }{\partial W_{k}}(\log (\det (\Gamma _{n}\left(W\right))))=\left(\left(\Gamma^{-1}_{n}\left(W\right)\right)_{ij}\right)_{1\leq i,j\leq d}^{T}\left(\frac{\left(\Gamma _{n}\left(W\right)\right)_{ij}}{\partial W_{k}}\right)_{1\leq i,j\leq d}\]
 with 
\[
\frac{\partial \Gamma _{ij}}{\partial W_{k}}=
\] 
\begin{equation}
\frac{1}{n}\sum _{t=1}^{n}\left[-\frac{\partial F_{W}(z_{t})(i)}{\partial W_{k}}\times \left(y_{t}-F_{W}(z_{t})\right)(j)-\frac{\partial F_{W}(z_{t})(j)}{\partial W_{k}}\left(y_{t}-F_{W}(z_{t})\right)(i)\right]\label{eq:deriv_{c}ov}
\end{equation}
so
\[
\frac{\partial }{\partial W_{k}}(\log (\det (\Gamma _{n}\left(W\right))))=\frac{1}{n}\left(\Gamma _{ij}^{-1}\right)_{1\leq i,j\leq d}^{T}\times
\]
\begin{equation}
\left(\sum _{t=1}^{n}-\frac{\partial F_{W}(z_{t})(i)}{\partial W_{k}}\times \left(y_{t}-F_{W}(z_{t})\right)(j)-\frac{\partial F_{W}(z_{t})(j)}{\partial W_{k}}\left(y_{t}-F_{W}(z_{t})\right)(i)\right)_{1\leq i,j\leq d}.\label{eq:sum_{p}ond}
\end{equation}

\begin{figure}[h]

\caption{\label{cap:MLP-restricted-to}MLP restricted to the output $2$ :
the plain lines}

\begin{center}\includegraphics[  scale=0.6]{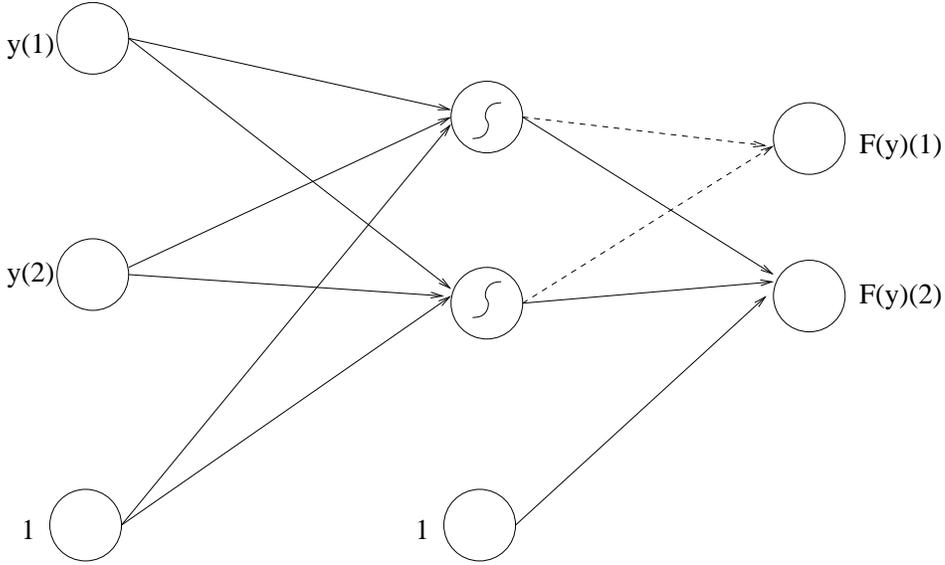}\end{center}
\end{figure}

The quantity $\frac{\partial F_{W}(z_{t})(i)}{\partial W_{k}}$ is
computed by back propagating the constant $1$ for the MLP
restricted to the output $i$. Figure \ref{cap:MLP-restricted-to}
gives an example of an MLP restricted to the output $2$.

Hence, the calculus of the gradient of $U_{n}\left(W\right)$ with respect
to the parameters of the MLP is straightforward. We have to compute
the derivative with respect to the weights of each single output MLP extracted from the original MLP by back propagating the constant value $1$,  then according to the
formula (\ref{eq:deriv_{c}ov}), we can compute easily the derivative
of each term of the empirical covariance matrix of the noise. Finally the gradient
is obtained by the sum of all the derivative terms of the empirical
covariance matrix multiplied by the terms of its inverse as in formula (\ref{eq:sum_{p}ond}).

\subsection{Second derivatives}
We write now 
\[
B_n(W_k,W_l):=\frac{1}{n}\sum _{t=1}^{n}\left( \frac{\partial F_{W}(z_{t})}{\partial W_{k}}\frac{\partial F_{W}(z_{t})}{\partial W_{l}}^T\right)
\]
and
\[
C_n(W_k,W_l):=\frac{1}{n}\sum _{t=1}^{n}\left( -(y_t-F_{W}(z_{t})) \frac{\partial^2 F_{W}(z_{t})}{\partial W_{k}\partial W_{l}}^T\right)
\]
We get

\[
\begin{array}{l}
\frac{\partial^2 U_n(W)}{\partial W_k\partial W_l}=\frac{\partial }{\partial W_{l}}2tr\left(\Gamma_{n}^{-1}(W)A_n(W_k)\right)=\\
2tr\left(\frac{\partial\Gamma_{n}^{-1}(W)}{\partial W_{l}}A(W_k)\right)+2tr\left(\Gamma_{n}^{-1}(W)B_n(W_k,W_l)\right)+2tr\left(\Gamma_{n}(W)^{-1}C_n(W_k,W_l)\right)\\
\end{array}
\]
Now, Magnus and Neudecker \cite{Magnus} give an analytic form of the derivative of an inverse matrix, from which we get
\[
\begin{array}{l}
\frac{\partial^2 U_n(W)}{\partial W_k\partial W_l}=2tr\left(\Gamma_{n}^{-1}(W)\left(A_n(W_k)+A_n^T(W_k)\right)\Gamma_{n}^{-1}(W)A_n(W_k)\right)+\\
2tr\left(\Gamma_{n}^{-1}(W)B_n(W_k,W_l)\right)+2tr\left(\Gamma_{n}^{-1}(W)C_n(W_k,W_l)\right)
\end{array}
\]
and
\begin{equation}\label{second_deriv}
\begin{array}{l}
\frac{\partial^2 U_n(W)}{\partial W_k\partial W_l}=4tr\left(\Gamma_{n}^{-1}(W)A_n(W_k)\Gamma_{n}^{-1}(W)A_n(W_k)\right)\\
+2tr\left(\Gamma_{n}^{-1}(W)B_n(W_k,W_l)\right)+2tr\left(\Gamma_{n}^{-1}(W)C_n(W_k,W_l)\right)
\end{array}
\end{equation}
\section{Asymptotic properties}
In the sequel, we will assume that the square of the noise $\varepsilon$ is integrable and that the cube of the variable $Z$ is integrable too. Moreover, it is easy to show that, for an MLP function, there exists a constant $C$ such that we have the following inequalities :
\[
\begin{array}{l}
\Vert \frac{\partial F_W(Z)}{\partial W_k}\Vert\leq C(1+\Vert Z\Vert)\\
\Vert \frac{\partial^2 F_W(Z)}{\partial W_k\partial W_l}\Vert\leq C(1+\Vert Z\Vert^2)\\
\Vert \frac{\partial^2 F_{W^1}(Z)}{\partial W_k\partial W_l}-\frac{\partial^2 F_{W^2}(Z)}{\partial W_k\partial W_l}\Vert\leq C\Vert W^1-W^2 \Vert(1+\Vert Z\Vert^3)
\end{array}
\]
These inequalities will be important to get the local asymptotic normality property implying the asymptotic normality of the parameter minimizing $U_n(W)$.
\subsection{Consistency and asymptotic normality of $\hat W_n$}
First we have to identify the contrast function associated with $U_n(W)$ 
\begin{Lem}\label{contraste}
\[
U_n(W)-U_n(W^0)\stackrel{a.s.}{\rightarrow}K(W,W^0)
\]
with $K(W,W^0)\geq 0$ and $K(W,W^0)= 0$ if and only if $W=W^0$.
\end{Lem}
\paragraph{Proof: }
Let us note 
\[
\Gamma(W)=E\left((Y-F_W(Z))(Y-F_W(Z))^T\right)
\]
the expectation of the covariance matrix of the noise for model parameter $W$.
By the strong law of large numbers we have
\[
\begin{array}{l}
U_n(W)-U_n(W^0)\stackrel{a.s.}{\rightarrow}\log \det(\Gamma(W))-\log \det(\Gamma(W^0))=\log\frac{\det(\Gamma(W))}{\det(\Gamma(W^0))}=\\
\log \det\left(\Gamma^{-1}(W^0) \left(\Gamma(W)-\Gamma(W^0)\right)+I_d\right)
\end{array}
\]
where $I_d$ denotes the identity matrix of $\mathbb R^d$.
So,  the lemma is true if $\Gamma(W)-\Gamma(W^0)$ is a positive matrix, null only if $W=W^0$. But this property is true since
\[
\begin{array}{l}
\Gamma(W)=E\left((Y-F_W(Z))(Y-F_W(Z))^T\right)=\\
E\left((Y-F_{W^0}(Z)+F_{W^0}(Z)-F_W(Z))(Y-F_{W^0}(Z)+F_{W^0}(Z)-F_W(Z))^T\right)=\\
E\left((Y-F_{W^0}(Z))(Y-F_{W^0}(Z))^T\right)+\\
E\left((F_{W^0}(Z)-F_W(Z))(F_{W^0}(Z)-F_W(Z))^T\right)=\\
\Gamma(W^0)+E\left((F_{W^0}(Z)-F_W(Z))(F_{W^0}(Z)-F_W(Z))^T\right)
\end{array}
\]
and the lemma follows from the identifiability assumption $\blacksquare$

We deduce  the theorem of consistency:
\begin{Theorem}
We have
\[
\hat W_n\stackrel{a.s.}{\rightarrow}W^0
\]
\end{Theorem} 
\paragraph{Proof} Remark that a constant $B$ exists such that 
\begin{equation}\label{borne}
sup_{W\in\Theta_s}\Vert Y-F_W(Z)\Vert^2<\Vert Y\Vert^2+B
\end{equation}
because $\Theta_s$ is compact, so $F_W(Z)$ is bounded. 
Let us define the function 
\[
\Phi(\Gamma):=\max(\log\det(\Gamma),d\log(\delta))
\]
where $d$ is the dimension of the observations $Y$ and $\delta>0$ strictly smaller than the smallest eigenvalue of $\Gamma_0$, since $\Gamma_0$ is definite positive we have for all $W$:
\[
\lim_{n\rightarrow\infty}\Phi(\Gamma_n(W))\stackrel{a.s.}{=}\lim_{n\rightarrow\infty}\log\det(\Gamma_n(W))=K(W,W^0)+\log\det(\Gamma^0)> d\log(\delta)
\] 
Now, for all $W$, thanks to the inequality (\ref{borne}) there exists constants $\alpha$ and $\beta$ such that 
\[
\left|\Phi\left((Y-F_W(Z))(Y-F_W(Z))^T\right)\right|\stackrel{a.s.}{<}\alpha\Vert Y\Vert^2+\beta
\]
but the right hand of this inequality is integrable, so the function $\Phi$ as an integrable envelope function
 and by  example 19.8 of van der Vaart \cite{Vandervaart} the set of functions $\left\{\Phi\left((Y-F_W(Z))(Y-F_W(Z))^T\right),\ W\in\Theta_s\right\}$ is Glivenko-Cantelli.

Now, the theorem 5.7 of van der Vaart \cite{Vandervaart}, shows  that $\hat W_n$ converges in probability to $W^0$, but it is easy to show that this convergence is almost sure. First, by lemma \ref{contraste}, we remark that for all neighborhood $\cal N$ of $W^0$ their exists a number $\eta({\cal N})>0$ such that for all $W\notin \cal N$ we have
\[
\log\det\left(\Gamma(W)\right)>\log\det\left(\Gamma(W^0)\right)+\eta({\cal N})
\]

   Now to show the strong consistency property we have to prove that for all neighborhood $\cal N$ of $W^0$ we have 
\(
\lim_{n\rightarrow \infty}\hat W_n \stackrel{a.s.}{\subset}{\cal N}
\)
or, equivalently, 
\[
\lim_{n\rightarrow \infty}\log\det\left(\Gamma(\hat W_n)\right)-\log\det\left(\Gamma(W^0)\right)<\eta({\cal N})
\]
By definition, we have 
\[
\log\det\left(\Gamma_n(\hat W_n)\right)\leq \log\det\left(\Gamma_n(W^0)\right)
\]
and the Glivenko-Cantelli property assures that
\[ 
\lim_{n\rightarrow \infty} \log\det\left(\Gamma_n(W^0)\right)-\log\det\left(\Gamma(W^0)\right)\stackrel{a.s.}{=} \lim_{n\rightarrow \infty} \Phi(\Gamma_n(W))-\log\det\left(\Gamma(W^0)\right)\stackrel{a.s.}{=}0
\]
therefore
\[
\lim_{n\rightarrow \infty} \log\det\left(\Gamma_n(\hat W_n)\right)<\log\det\left(\Gamma(W^0)\right)+\frac{\eta({\cal N})}{2}
\]
We have also
\[ 
\lim_{n\rightarrow \infty} \log\det\left(\Gamma_n(\hat W_n)\right)-\log\det\left(\Gamma(\hat W_n)\right)\stackrel{a.s.}{=} \lim_{n\rightarrow \infty} \Phi(\Gamma_n(W))-\log\det\left(\Gamma(\hat W_n)\right)\stackrel{a.s.}{=}0
\]
and finally
\[
\lim_{n\rightarrow \infty} \log\det\left(\Gamma(\hat W_n)\right)-\frac{\eta({\cal N})}{2}<
\log\det\left(\Gamma_n(\hat W_n)\right)<\log\det\left(\Gamma(W^0)\right)+\frac{\eta({\cal N})}{2}
\]
$\blacksquare$

Now, we can establish the asymptotic normality for the estimator. 
\begin{Lem}

Let $\Delta U_n(W^0)$ be the gradient vector of $U_n(W)$ at $W^0$,  $\Delta U(W^0)$ be the gradient vector of $U(W):=\log\det\left(\Gamma(W)\right)$ at $W^0$ and $HU_n(W^0)$ be the Hessian matrix of  $U_n(W)$ at $W^0$. 

We define finally
\[
B(W_k,W_l):=\frac{\partial F_{W}(Z)}{\partial W_{k}}\frac{\partial F_{W}(Z)}{\partial W_{l}}^T
\]

Then we get  
\begin{enumerate}
\item $HU_n(W^0)\stackrel{a.s.}{\rightarrow}2I_0$
\item $\sqrt{n}\Delta U_n(W^0)\stackrel{Law}{\rightarrow}{\cal N}(0,4I_0)$
\end{enumerate}
where, the component $(k,l)$ of the matrix $I_0$ is :
\[
tr\left(\Gamma^{-1}_0E\left(B(W^0_k,W^0_l)\right) \right)
\]
\end{Lem}
\paragraph{proof}
First we note
\[
A(W_k)=\left(-\frac{\partial F_{W}(Z)}{\partial W_{k}}(Y-F_{W}(Z))^T\right)
\]
To prove the lemma, we remark first that the component $(k,l)$ of the matrix $4I_0$ is :
\[
E\left(\frac{\partial U(W^0)}{\partial W_k}\frac{\partial U(W^0)}{\partial W^0_l}\right)=E\left(2tr\left(\Gamma^{-1}_0A^T(W^0_k)\right)\times2tr\left(\Gamma^{-1}_0A(W^0_l)\right)\right)
\]
and, since the trace of the product is invariant by circular permutation, 
\[
\begin{array}{l}
E\left(\frac{\partial U(W^0)}{\partial W_k}\frac{\partial U(W^0)}{\partial W^0_l}\right)=\\
4E\left( -\frac{\partial F_{W^0}(Z)^T}{\partial W_k}\Gamma^{-1}_0(Y-F_{W^0}(Z))(Y-F_{W^0}(Z))^T\Gamma^{-1}_0\left(-\frac{\partial F_{W^0}(Z))}{\partial W_l}\right)\right)\\
=4E\left(\frac{\partial F_{W^0}(Z)^T}{\partial W_k}\Gamma^{-1}_0\frac{\partial F_{W^0}(Z)}{\partial W_l}\right)\\
=4tr\left(\Gamma^{-1}_0E\left(\frac{\partial F_{W^0}(Z)}{\partial W_k}\frac{\partial F_{W^0}(Z)^T}{\partial W_l}\right) \right)\\
=4tr\left(\Gamma^{-1}_0E\left(B(W^0_k,W^0_l)\right)\right) 

\end{array}
\]
Now, for the component $(k,l)$ of the expectation of the Hessian matrix, we remark that 
\[
\lim_{n\rightarrow \infty}tr\left(\Gamma_{n}^{-1}(W^0)A_n(W^0_k)\Gamma_{n}^{-1}(W^0)A_n(W^0_k)\right)=0
\]
and
\[
\lim_{n\rightarrow \infty}tr\Gamma_{n}^{-1}C_n(W^0_k,W^0_l)=0
\]
so
\[
\begin{array}{l}
\lim_{n\rightarrow \infty}H_n(W^0)=\lim_{n\rightarrow \infty}4tr\left(\Gamma_{n}^{-1}(W^0)A_n(W^0_k)\Gamma_{n}^{-1}(W^0)A_n(W^0_k)\right)+\\
2tr\Gamma_{n}^{-1}(W^0)B_n(W^0_k,W^0_l)+2tr\Gamma_{n}^{-1}C_n(W^0_k,W^0_l)=\\
=2tr\left(\Gamma^{-1}_0E\left(B(W^0_k,W^0_l)\right)\right)\\
\blacksquare
\end{array}
\]

Now, from a classical argument of local asymptotic normality (see for example Yao \cite{Yao}), we deduce the following property for the estimator $\hat W_n$:  

\begin{Prop}
We have
\[
\lim_{n\rightarrow\infty}\sqrt{n}(\hat W_n-W^0)={\cal N}(0,I_0^{-1})
\]
\end{Prop}
However, if $W_n^*$ is the estimator of the generalized least squares :
\[
W_n^*:=\arg\min \frac{1}{n}\sum _{t=1}^{n}\left(Y_{t}-F_{W}\left(Z_{t}\right)\right)^{T}\Gamma_0^{-1}\left(Y_{t}-F_{W}\left(Z_{t}\right)\right)
\]
then we have also
\[
\lim_{n\rightarrow\infty}\sqrt{n}(W_n^*-W^0)={\cal N}(0,I_0^{-1})
\]
so $\hat W_n$ has the same asymptotic behavior as the generalized least squares estimator with the true covariance matrix  $\Gamma^{-1}_{0}$ which is asymptotically optimal (see for example Ljung \cite{Ljung}).  Therefore, the proposed estimator is asymptotically optimal too.
\subsection{Asymptotic distribution of the test statistic $T_n$}
Let us assume that the null hypothesis $H_0$ is true, we write\\  $\hat W_n=\arg\min_{W\in \Theta_s}U_n(W)$ and
$\hat W^0_n=\arg\min_{W\in \Theta_q} U_n(W)$, where $\Theta_q$ is viewed as a subset of $\Theta_s$. The asymptotic distribution of $T_n$ is then a consequence of the previous section. Namely, if we replace $nU_n(W)$ by its Taylor expansion around $\hat W_n$ and $\hat W^0_n$, following  van der Vaart \cite{Vandervaart} chapter 16 we have :
\[
T_n=\sqrt{n}\left(\hat W_n-\hat W^0_n\right)^T I_0\sqrt{n}\left(\hat W_n-\hat W^0_n\right)+o_P(1)\stackrel{\cal D}{\rightarrow}\chi^2_{s-q}
\]

\section{Experimental results}
\subsection{Simulated example}

Although the estimator associated with the cost function $U_{n}\left(W\right)$,
is theoretically better than the ordinary mean least squares estimator,
it is of some interest to quantify this fact by simulation. Moreover, there are some pitfalls in practical situations with MLPs. 

The first point is that we have no guaranty to reach the global minimum of the cost function, we can only hope to find a good local minimum if we are using many estimations with different initial weights.

The second point, is the fact that MLP are black box, it means that it is difficult to give an interpretation of their
parameters and it is almost impossible to compare MLP by comparing their parameters even if we try to take into account the possible permutations of the weights, because the difference between the weights may reflect only the differences of local minima reached by weights during the learning.

All these reasons explain why we choose, for simplicity,  to compare the estimated covariance
matrices of the noise instead of comparing directly the estimated parameters of MLPs.

\subsubsection{The model}

To simulate our data, we use an MLP with 2 inputs, 3 hidden units,
and 2 outputs. We choose to simulate an auto-regressive time series, where the outputs at time $t$ are the inputs
for time $t+1$. Moreover, with MLPs, the statistical properties
of such a model are the same as with independent identically distributed
(i.i.d.) data, because the time series constitutes a mixing process (see Yao \cite{Yao}). 

The equation of the model is the following \[
Y_{t+1}=F_{W_{0}}\left(Y_{t}\right)+\varepsilon _{t+1}\]

where 

\begin{itemize}
\item $Y_0=(0,0)$.
\item $\left(Y_{t}\right)_{1\leq t\leq 1000}$, $Y_{t}\in \mathbb{R}^{2}$,
is the bidimensional simulated random process 
\item $F_{W_{0}}$ is an MLP function with weights $W_{0}$ chosen randomly between $-2$ and 2.
\item $(\varepsilon _{t})$ is an i.i.d. centered noise with covariance matrix $\Gamma _{0}=\left(\begin{array}{cc}
 5 & 4\\
 4 & 5\end{array}
\right)$.
\end{itemize}
In order to study empirically the statistical properties of our estimator
we make $400$ independent simulations of the bidimensional time series
of length $1000$. 
\subsubsection{The results}
Our goal is to compare the estimator minizing $U_n\left(W\right)$ or equation (\ref{logdet}) and  the weights minimizing the mean square error (MSE), equation (\ref{ols}).
For each time series we estimate the weights of the MLP using the cost
function $U_n\left(W\right)$ and the MSE. 
The estimations have been done using the second order algorithm BFGS,
 and for each estimation we choose the best result obtained
after $100$ random initializations of the weights. Thus, we avoid
 plaguing our learning with poor local minima. 

We show here the mean of estimated covariance matrices of the noise for
$U_n(W)$ and the mean square error (MSE) cost function: \[
U_n\left(W\right)\, :\, \left(\begin{array}{cc}
 4.960 & 3.969\\
 3.969 & 4.962\end{array}
\right)\mbox {\, and\, MSE}\, :\, \left(\begin{array}{cc}
 4.938 & 3.932\\
 3.932 & 4.941\end{array}
\right)\]
The estimated standard deviation of the terms of the matrices are
all equal to $0.01$, so the differences observed between the two matrices are statistically significant. 
We can see that the estimated covariance of the noise is on average better with
the estimator associated to the cost function $U_n\left(W\right)$,
in particular it seems that there is slightly less overfitting with
this estimator, and the non diagonal terms are greater than with the least squares estimator. As expected, the determinant of the mean matrix associated with $U_n(W)$ is 8.86 instead of 8.93 for the matrix associated with the MSE. 
\subsection{Application to real time series: Pollution of ozone}
Ozone is a reactive oxidant, which is formed both in the stratosphere and troposphere. Near the ground's surface, ozone is directly harmful to human health, plant life and damages physical materials. The population, especially in large cities and in suburban zones which suffer from summer smog, wants to be warned of high pollutant concentrations in advance. The statistical ozone modelling and more particularly regression models have been widely studied \cite{Comrie}, \cite{Gardner}. Generally, linear models do not seem to capture all the complexity of the phenomena. Thus, the use of nonlinear techniques is recommended to deal with ozone prediction. Here we want to predict ozone pollution in two sites at the same time. The sites are the pollution levels in the south of Paris (13th district) and on the top of the Eiffel Tower. As these sites are very near each other we can expect that the two components of the noise are very correlated. 
\subsubsection{The model}
The neural model used in this study is autoregressive and includes exogenous parameters (called NARX model), where $X$ stands for exogeneous variables. Our aim is to predict the maximum level of ozone pollution of the next day knowing the today's maximum level of pollution and the maximal temperature of the next day. If we note $Y^1$ the level of pollution for Paris 13, $Y^2$ the level of pollution for the Eiffel Tower and $Temp$ the temperature, the model can be written as follows: 
\begin{equation}
(Y^1_{t+1},Y^2_{t+1})=F_W(Y^1_t,Y^2_t,Temp_{t+1})+\varepsilon_{t+1}
\end{equation}
We will assume that the variables are mixing as previously.
As usual with real time series, overtraining is a crucial problem. MLPs are very overparametrized models. This occurs when the model learns the details of the noise of the training data. Overtrained models have very poor performance on fresh data. To avoid overtraining we use in this study the SSM pruning technique,  a statistical stepwise method using a BIC-like criterion (Cottrell et al \cite{Cottrell}). The MLP with the minimal dimension is found by the elimination of the irrelevant weights. Here, we will compare behavior of this method for both cost function: The mean square error (MSE) and the logarithm of the determinant of the empirical covariance matrix of the noise ($U_n(W)$). 
\subsubsection{The dataset}
\begin{figure}[h]
\caption{\label{ozone}Ozone time series}
\begin{center}\includegraphics[scale=0.45,angle=-90]{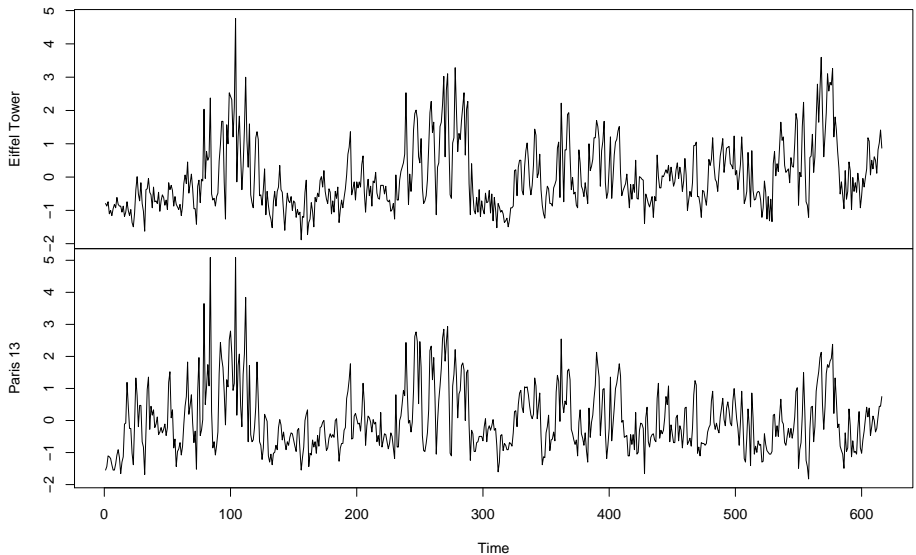}\end{center}
\end{figure}

This study presents the ozone concentration of the Air Quality Network of the Ile de France Region (AIRPARIF, Paris, France). The data used in this work are from 1994 to 1997, we use only the months from April to September inclusive because there is no peak during the winter period.  According to the model, we have the following parameters:
\begin{itemize}
\item The maximum temperature of the day
\item Persistence is used by introducing the previous day's peak ozone.
\end{itemize}
Before their use in the neural network, all these data have been centered and normalized. 
The data used to train the MLPs are chosen randomly in the whole period and we leave 100 observations  to form a fresh data set (test set), which will be used for models evaluation. In order to evaluate the models we repeat 400 times this random sampling to get 400 covariance matrices on each set for the two cost functions. 
Figure \ref{ozone} is a plot of the centered and normalized original data.

\subsubsection{The results}
\begin{figure}[h]
\caption{\label{predictions}Predicted time series}
\begin{center}\includegraphics[scale=0.45,angle=-90]{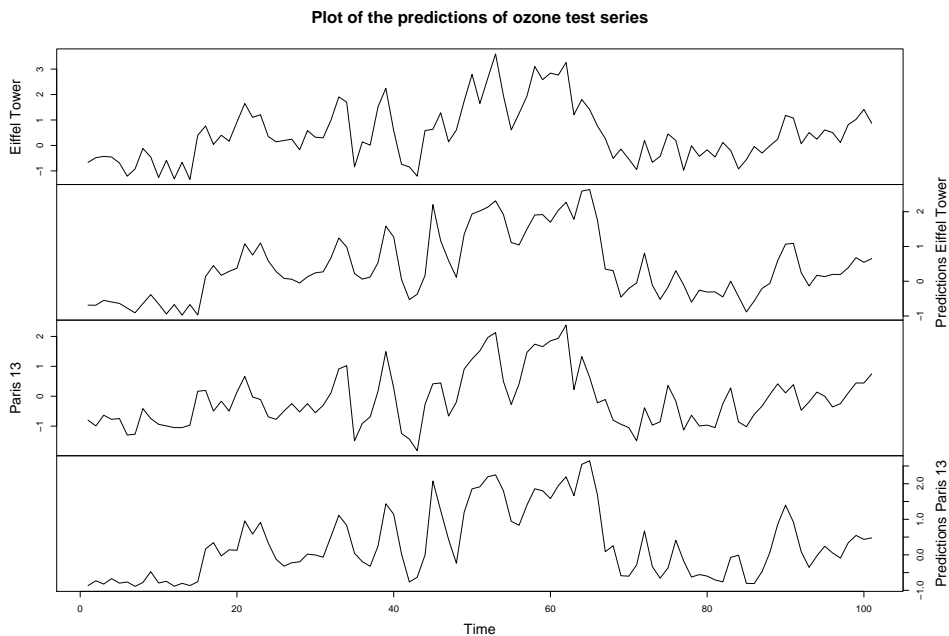}\end{center}
\end{figure}

For the learning set, we get the following results for the averaged covariance matrix (the estimated standard deviation for the coefficients is about 0.0005):
\[
U_n(W)\, :\, \left(\begin{array}{cc}
 0.27 & 0.20\\
 0.20 & 0.34\end{array}
\right)\mbox {\, and\, MSE}\, :\, \left(\begin{array}{cc}
 0.27 & 0.18\\
 0.18 & 0.34\end{array}
\right)\]
for the test set, we get the following results for the averaged covariance matrix (the estimated standard deviation for the coefficients is about 0.002) :
\[
U_n(W)\, :\, \left(\begin{array}{cc}
 0.29 & 0.22\\
 0.22 & 0.36\end{array}
\right)\mbox {\, and\, MSE}\, :\, \left(\begin{array}{cc}
 0.33 & 0.20\\
 0.20 &0.39\end{array}
\right)\]
The two matrices are almost the same for the learning set,  however the non-diagonal terms are greater for the $U_n(W)$ cost function. Moreover, looking at the averaged matrix on the test set, we see that the generalization capabilities  are better for $U_n(W)$ and the differences are statistically significant. 
Generally, the best MLP for $U_n(W)$ has less weights than the best MLP for the MSE cost function. Hence, the proposed cost function leads to a somewhat more parsimonious model, because the pruning technique is very sensitive to the variance of estimated parameters. This gain is valuable regarding the generalization capacity of the model, because the difference is almost null for the learning data set but is greater on the test data. Figure \ref{predictions} is a plot of the centered and normalized original test data and its prediction.
\section{Conclusion}
In the linear multidimensional regression model the optimal estimator has an analytic solution (see Magnus and Neudecker \cite{Magnus}), so it does not make sense to consider minimization of a cost function. However, for the non-linear multidimensional regression model, the ordinary least squares estimator is sub-optimal, if the covariance matrix of the noise is not the
identity matrix. We can overcome this difficulty by using the cost function $U_n(W)=\log\det(\Gamma_n(W))$.  In this paper, we have provided a proof of  the optimality of the estimator associated with  $U_n(W)$. Statistical thought tells us that it is always better for the neural networks practitionners to use a more efficient estimator because such estimator are better on average, even if the difference seems to be small. This estimator is especially important if the pratitionners are using pruning techniques. Indeed pruning technique are based on Wald test or approximated Wald test as for the optimal brain damage or optimal brain surgeon method (see Cottrell et al. \cite{Cottrell}) and these tests are very sensitive to the variance of the estimated parameters. Moreover, we have shown that this cost function leads to a simpler $\chi^2$ test to determine the number of weights if the model is identifiable. These theoretical results have been confirmed by a simulated example, and we have see for a real time series that we can expect slight improvement especially in model selection, this confirms the fact that such techniques are very sensitive to the variance of the estimated weights. 

\begin{footnotesize}

% IF YOU DO NOT USE BIBTEX, USE THE FOLLOWING SAMPLE SCHEME FOR THE REFERENCES
% ----------------------------------------------------------------------------

\end{footnotesize}

% ****************************************************************************
% END OF BIBLIOGRAPHY AREA
% ****************************************************************************


\begin{thebibliography}{99}

\bibitem{Comrie}
A.C. Comrie, 
Comparing neural networks and regression models for ozone forecasting,
Air and Waste Management Association, 47 (1997) 653-663.

\bibitem{Cottrell}
M. Cottrell,et al.,
Neural modeling for time series: a statistical stepwise method for weight elimination,
IEEE Transaction on Neural Networks 6 (1995) 1355-1364.

\bibitem {Fukumizu}
K. Fukumizu, 
A regularity condition of the information matrix of a multilayer perceptron network,
Neural Networks, 9:5 (1996) 871-879. 

\bibitem {Fukumizu2}
K. Fukumizu, 
Likelihood ratio of unidentifiable models and multilayer neural networks,
The Annals of Statistics, 31:3 (2003) 833-851. 


\bibitem {Gallant}
R. A. Gallant,
Non linear statistical models (J. Wiley and Sons, New York,1987). 

\bibitem{Gardner}
M.W. Gardner and S.R. Dorling, 
Artificial neural networks, the multilayer Perceptron. A review of applications in the atmospheric sciences,
Atmospheric Environment, 32:14/15 (1998) 2627-2636.

\bibitem {Ljung}
L. Ljung,
System identification: Theory for the user
(Prentice Hall, New Jersey, 1999).

\bibitem {Magnus}
J. Magnus and H. Neudecker,
Matrix differential calculus with applications in statistics and econometrics
(J. Wiley and Sons, New York, 1988).

\bibitem{Rynkiewicz4}
 J. Rynkiewicz,
 Estimation of Multidimensional Regression Model with Multilayer Perceptron, in: J. Mira and A. Prieto, ed., proc. IWANN'2003, Lecture Notes in Computer Science, Vol. 2686 (Springer, Berlin, 2003) 310-317.

\bibitem{Sussmann}
H.J. Sussmann,
Uniqueness of the weights for minimal feedforward nets with a given input-output Map,
Neural Networks 5 (1992) 589-593.

\bibitem{Vandervaart}
A. W. Van der Vaart,
Asymptotic statistics
(Cambridge University Press, Cambridge, 1998).

\bibitem{Williams}
P. M. Williams,
Using neural networks to model conditional multivariate densities,
Neural Computation 8:4 (1996) 843-854.

\bibitem{Yao}
J.F. Yao,
On least squares estimation for stable nonlinear AR processes,
The Annals of the Institute of Mathematical Statistics 52 (2000) 316-331.  
\end{thebibliography}
\end{document}